\documentclass[15pt, a4paper]{article}
\usepackage{}
\usepackage{amsfonts}
\usepackage{mathbbold}
\usepackage{bbm}
\usepackage{mathrsfs}
\usepackage{latexsym}
\usepackage{graphicx}
\usepackage{amssymb}

\newtheorem{theorem}{Theorem}[section]
\newtheorem{lemma}[theorem]{Lemma}

\title{{\Large \bf The unicyclic hypergraph with extremal spectral radius\thanks{Supported by National natural science foundation of China (NSFC)
(Nos. Nos. 12171222, 12101285, 12071411), Natural science foundation of Guangdong province (No. 2021A1515010254), Foundation of Lingnan Normal
University(ZL1923).}}}
\author{Guanglong Yu$^{a}$
~ Lin Sun$^{a}$\thanks{Corresponding authors, E-mail addresses:
yglong01@163.com (G. Yu), fiona$_{-}$sl@163.com (L. Sun).} ~  Hailiang Zhang$^b$ ~ Gang Li$^{c}$ ~
\\ ~ \\
{\footnotesize $^a$Department of Mathematics, Lingnan Normal
University,  Zhanjiang, Guangdong, 524048, P.R.China}\\
{\footnotesize $^b$Department of Mathematics, Taizhou University, Linhai, Zhejiang, 317000, P.R.China}\\
{\footnotesize $^c$College of Mathematics and Systems Science, Shandong University of Science and Technology,}\\
{\footnotesize  Qingdao, Shandong, 266590, P.R.China}}
\date{}

\begin{document}
\maketitle

\begin{abstract}
For a $hypergraph$ $\mathcal{G}=(V, E)$ consisting of a nonempty vertex set $V=V(\mathcal{G})$ and an edge set $E=E(\mathcal{G})$, its $adjacency$ $matrix$ $\mathcal {A}_{\mathcal{G}}=[(\mathcal {A}_{\mathcal{G}})_{ij}]$ is defined as
$(\mathcal {A}_{\mathcal{G}})_{ij}=\sum_{e\in E_{ij}}\frac{1}{|e| - 1}$, where $E_{ij} = \{e \in E \, |\,  i, j \in e\}$.
The $spectral$ $radius$ of a hypergraph $\mathcal{G}$, denoted by $\rho(\mathcal {G})$, is the maximum modulus among all eigenvalues of $\mathcal {A}_{\mathcal{G}}$.
In this paper, among all $k$-uniform ($k\geq 3$) unicyclic hypergraphs with fixed number of vertices, the hypergraphs with the maximum and the second the maximum spectral radius are completely determined, respectively.


\end{abstract}
\baselineskip 18.6pt

\section{Introduction}

\ \ \ \
In the past twenty years, different hypermatrices or tensors for hypergraphs have been developed to explore spectral hypergraph theory. Many interesting spectral properties of hypergraphs have been explored \cite{Ban.CB}-\cite{COQY}, \cite{LQIE}-\cite{JZLS}.
Recently, A. Banerjee \cite{A.Ban} introduced an adjacency matrix and use its spectrum so that some spectral and structural
properties of hypergraphs are revealed. In this paper, we go on studying the spectra of hypergraphs according to the adjacency matrix introduced in \cite{A.Ban}.

Now we recall some notations and definitions related to hypergraphs.
For a set $S$, we denote by $|S|$ its cardinality. A $hypergraph$ $\mathcal{G}=(V, E)$ consists of a nonempty vertex set $V=V(\mathcal{G})$ and an edge set $E=E(\mathcal{G})$, where each edge $e\in E(\mathcal{G})$ is a subset of $V(\mathcal{G})$ containing at least two vertices. The cardinality $n=|V(\mathcal{G})|$ is called the order; $m=|E(\mathcal{G})|$ is called the edge number of hypergraph $\mathcal{G}$. Denote by $t$-set a set with size (cardinality) $t$. We say that a hypergraph $\mathcal{G}$ is $uniform$ if its every edge has the same size, and call it $k$-$uniform$ if its every edge has size $k$ (i.e. every edge is a $k$-subset). It is known that a $2$-uniform graph is always called a ordinary graph or graph for short.

For a hypergraph $\mathcal{G}$, we define $\mathcal{G}-e$ ($\mathcal{G}+e$)
to be the hypergraph obtained from $\mathcal{G}$ by deleting the edge $e\in
 E(\mathcal{G})$ (by adding a new edge $e$ if $e\notin
 E(\mathcal{G})$); for an edge subset $B\subseteq E(\mathcal{G})$, we define $\mathcal{G}-B$
to be the hypergraph obtained from $\mathcal{G}$ by deleting each edge $e\in
 B$; for a vertex subset $S\subseteq V(\mathcal{G})$, we define $\mathcal{G}-S$ to be the hypergraph obtained from $\mathcal{G}$ by deleting all the vertices in $S$ and deleting the edges incident with any vertex in $S$. For two $k$-uniform hypergraphs $\mathcal{G}_{1}=(V_{1}, E_{1})$ and $\mathcal{G}_{2}=(V_{2}, E_{2})$, we say the two graphs are $isomorphic$ if there is a bijection $f$ from $V_{1}$ to $V_{2}$, and there is a bijection $g$ from $E_{1}$ to $E_{2}$ that maps each edge $\{v_{1}$, $v_{2}$, $\ldots$, $v_{k}\}$ to $\{f(v_{1})$, $f(v_{2})$, $\ldots$, $f(v_{k})\}$.

In a hypergraph, two vertices are said to be $adjacent$ if
both of them are contained in an edge. Two edges are said to be $adjacent$
if their intersection is not empty. An edge $e$ is said to be $incident$ with a vertex $v$ if
$v\in e$. The $neighbor$ $set$ of vertex $v$ in hypergraph $\mathcal{G}$, denoted by $N_{\mathcal{G}}(v)$, is the set of vertices adjacent to $v$ in $\mathcal{G}$. The $degree$ of a vertex $v$ in $\mathcal{G}$, denoted by $deg_{\mathcal{G}}(v)$ (or $deg(v)$ for short), is the number of the edges incident with $v$. A vertex of degree $1$ is called a $pendant$ $vertex$. A $pendant$ $edge$ is an edge with at most one vertex of degree more than one and other vertices in this edge being all pendant vertices.

In a hypergraph, a $hyperpath$ of length $q$ ($q$-$hyperpath$) is defined to be an alternating sequence
of vertices and edges $v_{1}e_{1}v_{2}e_{2}\cdots v_{q}e_{q}v_{q+1}$ such that
(1) $v_{1}$, $v_{2}$, $\ldots$, $v_{q+1}$ are all distinct vertices;
(2) $e_{1}$, $e_{2}$, $\ldots$, $e_{q}$ are all distinct edges;
(3) $v_{i}$, $v_{i+1}\in e_{i}$ for $i = 1$, $2$, $\ldots$, $q$;
(4) $e_{i}\cap e_{i+1}=v_{i+1}$ for $i = 1$, $2$, $\ldots$, $q-1$; (5) $e_{i}\cap e_{j}=\emptyset$ if $|i-j|\geq 2$.
If there is no discrimination, a hyperpath is sometimes written as $e_{1}e_{2}\cdots e_{q-1}e_{q}$, $e_{1}v_{2}e_{2}\cdots v_{q}e_{q}$ or $v_{1}e_{1}v_{2}e_{2}\cdots v_{q}e_{q}$. A $hypercycle$ of length $q$ ($q$-$hypercycle$) $v_{1}e_{1}v_{2}e_{2}\cdots v_{q-1}e_{q-1}v_{q}e_{q}v_{1}$ is obtained from a hyperpath $v_{1}e_{1}v_{2}e_{2}\cdots v_{q-1}e_{q-1}v_{q}$ by adding a new edge $e_{q}$ between $v_{1}$ and $v_{q}$ where $e_{q}\cap e_{1}= \{v_{1}\}$, $e_{q}\cap e_{q-1}= \{v_{q}\}$, $e_{q}\cap e_{j}=\emptyset$ if $j\neq1, q-1$ and $|q-j|\geq 2$. The length of a hyperpath $P$ (or a hypercycle $C$), denoted by $L(P)$ (or $L(C)$), is the number of the edges in $P$ (or $C$). In a hypergraph $\mathcal{G}$, the $distance$ between vertices $u$ and $v$, denoted by $dist_{\mathcal{G}}(u,v)$, is the length of the shortest hyperpath between $u$ and $v$. Denote by $\mathcal{C}_{n,k}$ a $k$-uniform cycle of order $n$. A hypergraph $\mathcal{G}$
is connected if there exists a hyperpath from $v$ to $u$ for all $v, u \in V$,
and $\mathcal{G}$ is called $acyclic$ if it contains no hypercycle.

Recall that a tree is an ordinary graph which is 2-uniform, connected and acyclic. A $supertree$ is similarly defined to be a hypergraph which is both connected and acyclic, which is a natural generalization of the ordinary tree.
Clearly, in a supertree, its each pair of the
edges have at most one common vertex. Therefore, the edge number of a $k$-uniform supertree of order $n$ is $m=\frac{n-1}{k-1}$.

A connected $k$-uniform hypergraph with $n$ vertices and $m$ edges is $r$-$cyclic$ if $n-1 =(k-1)m-r$. For $r=1$, it is called a $k$-uniform $unicyclic$ $hypergraph$; for $r=0$, it is a $k$-uniform supertree.
In \cite{COQY}, C. Ouyang et al. proved that a simple connected $k$-graph $G$ is unicyclic (1-cyclic) if and only if it has only one cycle. From this, for a unicyclic hypergraph $G$ with unique cycle $C$, it follows that (1) if $L(C)=2$, then the two edges in $C$ have exactly two common vertices, and $|e\cap f|\leq 1$ for any two edges $e$ and $f$ not in $C$ simultaneously; (2) if $L(C)\geq 3$, then any two edges in $G$ have at most one common vertices; (3) every connected component of $G-E(C)$ is a supertree.

\setlength{\unitlength}{0.65pt}
\begin{center}
\begin{picture}(575,250)
\put(370,118){\circle*{4}}
\put(504,118){\circle*{4}}
\qbezier(370,118)(440,250)(504,118)
\qbezier(370,118)(439,173)(504,118)
\qbezier(370,118)(436,6)(504,118)
\qbezier(370,118)(440,73)(504,118)
\qbezier(504,118)(484,182)(512,181)
\qbezier(512,181)(542,173)(504,118)
\qbezier(504,118)(558,161)(570,139)
\qbezier(504,118)(575,108)(570,139)
\qbezier(504,118)(483,59)(508,60)
\qbezier(504,118)(540,65)(508,60)
\qbezier(370,118)(300,149)(298,118)
\qbezier(298,118)(301,88)(370,118)
\put(21,110){\circle*{4}}
\put(161,111){\circle*{4}}
\qbezier(21,110)(94,232)(161,111)
\qbezier(21,110)(91,150)(161,111)
\qbezier(21,110)(94,70)(161,111)
\qbezier(21,110)(94,0)(161,111)
\put(199,88){\circle*{4}}
\put(203,95){\circle*{4}}
\put(195,80){\circle*{4}}
\put(543,95){\circle*{4}}
\put(538,86){\circle*{4}}
\put(548,104){\circle*{4}}
\put(88,26){$\mathscr{U}^{\ast}$}
\put(433,26){$\mathcal {F}$}
\put(235,-6){Fig. 1.1. $\mathscr{U}^{\ast}, \mathcal {F}$}
\put(53,139){$e_{2}$}
\put(3,110){$v_{2}$}
\put(53,81){$e_{1}$}
\put(171,146){$e_{3}$}
\put(195,117){$e_{4}$}
\put(137,109){$v_{1}$}
\put(358,125){$v_{2}$}
\put(307,115){$e'_{3}$}
\put(481,116){$v_{1}$}
\put(406,83){$e_{1}$}
\put(400,151){$e_{2}$}
\put(504,165){$e_{4}$}
\qbezier(161,111)(159,168)(183,165)
\qbezier(161,111)(212,152)(183,165)
\qbezier(161,111)(214,144)(221,123)
\qbezier(161,111)(228,94)(221,123)
\qbezier(161,111)(154,53)(180,59)
\qbezier(161,111)(205,72)(180,59)
\put(168,71){$e_{m}$}
\put(546,129){$e_{5}$}
\put(499,76){$e_{m}$}
\end{picture}
\end{center}

Let $\mathbb{C}=v_{1}e_{1}v_{2}e_{2}v_{1}$ be the $k$-uniform 2-hypercycle with $e_{i}=\{v_{1}$, $v_{a(i,1)}$, $v_{a(i,2)}$, $\ldots$, $v_{a(i,k-2)}$, $v_{2}\}$ for $i=1, 2$. Let $\mathscr{U}^{\ast}$ be the unicyclic hypergraph  of order $n$ obtained from cycle $\mathbb{C}$ by attaching $\frac{n}{k-1}-2$ pendant edges at vertex $v_{1}$. $\mathcal {F}$ be the $k$-uniform unicyclic hypergraph of order $n$ consisting of hypercycle $\mathbb{C}$ and $\frac{n}{k-1}-3$ pendant edges attaching at vertex $v_{1}$ and one pendant edge attaching at vertex $v_{2}$.

Let $E_{ij} = \{e \in E \, |\, i, j \in e\}$. The $adjacency$ $matrix$ $\mathcal {A}_{\mathcal{G}}=[(\mathcal {A}_{\mathcal{G}})_{ij}]$ of a hypergraph $\mathcal{G}$ is defined as
$$(\mathcal {A}_{\mathcal{G}})_{ij}=\sum_{e\in E_{ij}}\frac{1}{|e| - 1}.$$ It is easy to find that $\mathcal {A}_{\mathcal{G}}$ is symmetric if there is no requirement for direction on hypergraph $\mathcal{G}$, and find that $\mathcal {A}_{\mathcal{G}}$ is very convenient to be used to investigate the spectum of a hypergraph even without the requirement for edge uniformity. The $spectral$ $radius$ $\rho(\mathcal {G})$ of a hypergraph $\mathcal{G}$ is defined to be the spectral radius $\rho(\mathcal {A}_{\mathcal{G}})$, which is the maximum modulus among all eigenvalues of $\mathcal {A}_{\mathcal{G}}$.
In spectral theory of hypergraphs, the spectral radius is an index that attracts much attention due to its fine properties \cite{KCKPZ, YFTP, SFSH, LSQ, HLBZ, LLSM, COQY, YQY}.

We assume that the hypergraphs throughout
this paper are simple, i.e. $e_{i} \neq e_{j}$ if $i \neq j$, and assume the hypergraphs throughout
this paper are undirected. In this paper, among all $k$-uniform ($k\geq 3$) unicyclic hypergraphs with fixed number of vertices, the hypergraphs with the maximum and the second the maximum spectral radius are completely determined, respectively, getting the following result:

\begin{theorem}\label{th01.01} 
Let $\mathcal{G}$ be a $k$-uniform ($k\geq 3$) unicyclic hypergraph of order $n$. Then $\rho(\mathcal{G})\leq\rho(\mathscr{U}^{\ast})$ with equality if and only if $\mathcal{G}\cong \mathscr{U}^{\ast}$.

\end{theorem}

Let $\Lambda=\{\mathcal{G}\, |\, \mathcal{G}$ be a $k$-uniform unicyclic  ($k\geq 3$) hypergraph of order $n$, $m(\mathcal{G})\geq 4$ and $\mathcal{G}\ncong \mathscr{U}^{\ast}\}$.

\begin{theorem}\label{th01.02} 
Let $\mathcal{G}\in \Lambda$. Then
$\rho(\mathcal{G})\leq\rho(\mathcal {F})$ with equality if and only if $\mathcal{G}\cong \mathcal {F}$.

\end{theorem}

The layout of this paper is as follows: section 2 introduces some basic knowledge; section 3 represents our results.

\section{Preliminary}

\ \ \ \ \ For the requirements in the narrations afterward, we need some prepares. For a hypergraph $\mathcal{G}$ with vertex set $\{v_{1}$, $v_{2}$, $\ldots$, $v_{n}\}$, a vector on $\mathcal{G}$ is a vector $X=(x_{v_1}, x_{v_2}, \ldots, x_{v_n})^T \in R^n$ that
 entry $x_{v_i}$ is mapped to vertex $v_i$ for $i\leq i\leq n$.

 From \cite{OP}, by the famous Perron-Frobenius theorem, for $\mathcal {A}_{\mathcal{G}}$ of a connected uniform hypergraph $\mathcal{G}$ of order $n$, we know that there is unique one positive eigenvector $X=(x_{v_{1}}$, $x_{v_{2}}$, $\ldots$, $x_{v_{n}})^T \in R^{n}_{++}$ ($R^{n}_{++}$ means the set of positive real vectors of dimension $n$) corresponding to $\rho(\mathcal{G})$, where $\sum^{n}_{i=1}x^{2}_{v_{i}}= 1$ and each entry $x_{v_i}$ is mapped to each vertex $v_i$ for $i\leq i\leq n$. We call such an eigenvector $X$
the $principal$ $eigenvector$ of $\mathcal{G}$.

Let $A$ be an irreducible nonnegative $n \times n$ real matrix (with every entry being real number) with spectral radius $\rho(A)$. The following extremal representation (Rayleigh quotient) will be useful:
$$\rho(A)=\max_{X\in R^{n}, X\neq0}\frac{X^{T}AX}{X^{T}X},$$ and if a vector $X$ satisfies that $\frac{X^{T}AX}{X^{T}X}=\rho(A)$, then $AX=\rho(A)X$.

\section{Main results}

~~~~Recall that an automorphism of a $k$-uniform hypergraph $G$ is a permutation $\sigma$ of
$V(G)$ such that $\{i_{1}$, $i_{2}$, $\ldots$, $i_{k}\} \in E(G)$ if and only if $\{\sigma(i_{1}), \sigma(i_{2}), \ldots, \sigma(i_{k})\} \in
E(G)$, for any $i_{j}\in V(G)$, $j = 1$, $2$, $\ldots$, $k$. We denote by $Aut(G)$ the group of all automorphisms of $G$.

Denote by $P_{\sigma}$ a permutation matrix of order $n$ corresponding to a permutation $\sigma\in S_{n}$, where $S_{n}$ is the sysmmetric group of degree n. It is known that, for two matrices $\mathcal {D}, \mathcal {A}$ of order $n$ that $\mathcal {D} = P_{\sigma}\mathcal {A}P_{\sigma}$ for some permutation matrix $P_{\sigma}$, then $\mathcal {D}_{ij} = \mathcal {A}_{\sigma(i)\sigma(j)}$ in
$\mathcal {D}$.

\begin{lemma}\label{le03,00} 
Let $G$ be a $k$-uniform hypergraph on $n$ vertices. A permutation $\sigma \in S_{n}$ is an automorphism of
$G$ if and only if $P_{\sigma}\mathcal {A}_{G} =\mathcal {A}_{G}P_{\sigma}$.
\end{lemma}

\begin{proof}
Note that $P_{\sigma}$ is a permutation matrix. Then $P_{\sigma}^{2}=I$ (the identity matrix). Let $B=P_{\sigma}\mathcal {A}_{G}P_{\sigma}$.
Permutation $\sigma \in S_{n}$ is an automorphism of
$G$ if and only if $B=\mathcal {A}_{G}$. Thus we have $P_{\sigma}\mathcal {A}_{G} =\mathcal {A}_{G}P_{\sigma}$.
This completes the proof. \ \ \ \ \ $\Box$
\end{proof}

\begin{lemma}\label{le03,01} 
Let $G$ be a connected $k$-uniform hypergraph, $\mathcal {A}_{G}$
be its
adjacency matrix. If $\mathcal {A}X=\lambda X$, then for each automorphism $\sigma$ of $G$, we have $\mathcal {A} (P_{\sigma} X)=\lambda (P_{\sigma}X)$.
Moreover, if $\mathcal {A}X=\rho(\mathcal {A}) X$, then

$\mathrm{(1)}$ $P_{\sigma} X = X$;

$\mathrm{(2)}$ for any orbit $\Omega$ of $Aut (G)$ and each pair of vertices $i, j\in \Omega$, $x_{i}=x_{j}$ in $X$.
\end{lemma}

\begin{proof}
If $\mathcal {A}X=\lambda X$, then for each automorphism $\sigma$ of $G$, using Lemma \ref{le03,00}
gets $$\mathcal {A} (P_{\sigma} X)=P_{\sigma}\mathcal {A} X=P_{\sigma}(\lambda X)=\lambda (P_{\sigma}X).$$

If $\mathcal {A}X=\rho(\mathcal {A}) X$, since
$\mathcal {A}_{G}$ is nonnegative irreducible, then by Perron-Frobenius theorem \cite{OP}, we have $P_{\sigma} X = X$, and for any orbit $\Omega$ of $Aut (G)$ and each pair of vertices $i, j\in \Omega$, we have $x_{i}=x_{j}$ in $X$.
This completes the proof. \ \ \ \ \ $\Box$
\end{proof}

\begin{lemma}{\bf \cite{GYSLZ}}\label{le03.02} 
Let $e_{1}=\{v_{1,1}$, $v_{1,2}$, $\ldots$, $v_{1,k}\}$, $e_{2}=\{v_{2,1}$, $v_{2,2}$, $\ldots$, $v_{2,k}\}$, $\ldots$, $e_{j}=\{v_{j,1}$, $v_{j,2}$, $\ldots$, $v_{j,k}\}$ be some edges in a connected $k$-uniform hypergraph $\mathcal{G}$; $v_{u,1}$, $v_{u,2}$, $\ldots$, $v_{u,t}$ be vertices in $\mathcal{G}$ that $t< k$. For $1\leq i\leq j$, $\{v_{u,1}, v_{u,2}, \ldots, v_{u,t}\}\nsubseteq e_{i}$, $e^{'}_{i}=(e_{i}\setminus \{v_{i,1}, v_{i,2}, \ldots, v_{i,t}\})\cup\{v_{u,1}, v_{u,2}, \ldots, v_{u,t}\}$ satisfying that $e^{'}_{i}\notin E(\mathcal{G})$. Let $\mathcal{G}^{'}=\mathcal{G}-\sum e_{i}+\sum e^{'}_{i}$. If in the principal eigenvector $X$ of $\mathcal{G}$, for $1\leq i\leq j$, $x_{v_{i,1}}\leq x_{v_{u,1}}$, $x_{v_{i,2}}\leq x_{v_{u,2}}$, $\ldots$, and $x_{v_{i,t}}\leq x_{v_{u,t}}$, then $\rho(\mathcal{G}^{'})>\rho(\mathcal{G})$.
\end{lemma}

\begin{Proof}
Suppose $\mathcal{U}$ is a $k$-uniform $(k\geq 3)$ uicyclic hypergraph of order $n$ satisfying that $\rho(\mathcal{U})=\max\{\rho(\mathcal{G})\, |\ \mathcal{G}$ is a $k$-uniform $(k\geq 3)$ uicyclic hypergraph of order $n\}$.
Let $X$ be
the principal eigenvector of $\mathcal{U}$ and $x_{u}=\max\{x_{v}\, |\, v\in V(\mathcal{U})\}$. Denote by $\mathcal{C}=v_{1}e_{1}v_{2}e_{2}\cdots v_{q-1}e_{q-1}v_{q}e_{q}v_{1}$ the unique hypercycle in $\mathcal{U}$ with $e_{i}=\{v_{i}$, $v_{a(i,1)}$, $v_{a(i,2)}$, $\ldots$, $v_{a(i,k-2)}$, $v_{i+1}\}$ for $1\leq i\leq q-1$, $e_{q}=\{v_{q}$, $v_{a(q,1)}$, $v_{a(q,2)}$, $\ldots$, $v_{a(q,k-2)}$, $v_{1}\}$.

{\bf Assertion 1} $q=2$. Otherwise, suppose $q>2$. We suppose that $x_{\omega}=\max\{x_{v_{i}}\, |\ 1\leq i\leq q, v_{i}\in V(\mathcal{C})\}$. Without loss of generality, suppose $\omega=q$. Now, let $e^{'}_{1}=\{v_{1}$, $v_{a(1,1)}$, $v_{a(1,2)}$, $\ldots$, $v_{a(1,k-2)}$, $v_{q}\}$, $\mathcal{U}^{'}=\mathcal{U}-e_{1}+e^{'}_{1}$. Note that $\mathcal{U}^{'}$ is also a unicyclic hypergraph. Using Lemma \ref{le03.02} gets that $\rho(\mathcal{U}^{'})>\rho(\mathcal{U})$, which contradicts the maximality of $\rho(\mathcal{U})$.

{\bf Assertion 2} There exists $x_{\mu}=x_{u}$ where $\mu\in V(\mathcal{C})$. Otherwise, suppose $x_{v}<x_{u}$ for any $v\in V(\mathcal{C})$. Suppose $ue_{s_{1}}v_{s_{1}}e_{s_{2}}\cdots e_{s_{\varepsilon}}v_{s_{\varepsilon}}$ is the hyperpath from $u$ to $\mathcal{C}$ where $v_{s_{\varepsilon}}\in V(\mathcal{C})$. Note that $q=2$. Next, we consider two cases.

{\bf Case 1} $v_{s_{\varepsilon}}\in \{v_{1}, v_{2}\}$. Without loss of generality, suppose $v_{s_{\varepsilon}}=v_{1}$. Let $e^{'}_{1}=\{u$, $v_{a(1,1)}$, $v_{a(1,2)}$, $\ldots$, $v_{a(1,k-2)}$, $v_{2}\}$, $e^{'}_{2}=\{v_{2}$, $v_{a(2,1)}$, $v_{a(2,2)}$, $\ldots$, $v_{a(2,k-2)}$, $u\}$, $\mathcal{U}^{'}=\mathcal{U}-e_{1}+e^{'}_{1}-e_{2}+e^{'}_{2}$. As Assertion 1, it proved that $\rho(\mathcal{U}^{'})>\rho(\mathcal{U})$, which contradicts the maximality of $\rho(\mathcal{U})$.

{\bf Case 2} $v_{s_{\varepsilon}}\in (V(\mathcal{C})\setminus\{v_{1}, v_{2}\})$.  Without loss of generality, suppose $v_{s_{\varepsilon}}=v_{a(1,1)}$. Let $e^{'}_{1}=\{v_{1}$, $u$, $v_{a(1,2)}$, $\ldots$, $v_{a(1,k-2)}$, $v_{2}\}$, $\mathcal{U}^{'}=\mathcal{U}-e_{1}+e^{'}_{1}$. As Assertion 1, it proved that $\rho(\mathcal{U}^{'})>\rho(\mathcal{U})$, which contradicts the maximality of $\rho(\mathcal{U})$.

From the above two cases, then Assertion 2 follows as desired. By Assertion 2, we suppose $x_{\eta}=x_{u}$ where $\eta\in V(\mathcal{C})$. Then we have the following Claim 1.

{\bf Claim 1} $\eta\in \{v_{1}, v_{2}\}$. Otherwise, suppose $\eta\notin \{v_{1}, v_{2}\}$, and without loss of generality, we suppose that $\eta=v_{a(1,1)}$. Let $e^{'}_{2}=\{v_{2}$, $v_{a(2,1)}$, $v_{a(2,2)}$, $\ldots$, $v_{a(2,k-2)}$, $\eta\}$, $\mathcal{U}^{'}=\mathcal{U}-e_{2}+e^{'}_{2}$. As Assertion 1, it proved that $\rho(\mathcal{U}^{'})>\rho(\mathcal{U})$, which contradicts the maximality of $\rho(\mathcal{U})$. Thus, our Claim 1 holds.

 By Claim 1 above, without loss of generality, we suppose $x_{v_{1}}=x_{u}$. Then we have the following Claim 2.

{\bf Claim 2} $dist_{\mathcal{U}}(v_{1}, v)=1$ for any $v\in (V(\mathcal{U})\setminus \{v_{1}\})$. Otherwise, suppose $v_{\varphi}\in (V(\mathcal{U})\setminus \{v_{1}\})$, but $dist(v_{\varphi}, v_{1})\geq 2$. Obviously, $v_{\varphi}\notin V(\mathcal{C})$. Note that every connected component of $G-E(\mathcal{C})$ is a supertree. Suppose $v_{\varphi}e_{t_{1}}v_{t_{1}}e_{t_{2}}\cdots e_{t_{\vartheta}}v_{t_{\vartheta}}$ is the hyperpath from $v_{\varphi}$ to $\mathcal{C}$ where $v_{t_{\vartheta}}\in V(\mathcal{C})$. Next, we consider two cases.

{\bf Case 1} $v_{t_{\vartheta}}=v_{1}$. Then $\vartheta\geq 2$. Denote by $e_{t_{\vartheta-1}}=\{v_{t_{\vartheta-2}}$, $v_{a(t_{\vartheta-1},1)}$, $v_{a(t_{\vartheta-1},2)}$, $\ldots$, $v_{a(t_{\vartheta-1},k-2)}$, $v_{t_{\vartheta-1}}\}$. Let $e^{'}_{t_{\vartheta-1}}=\{v_{t_{\vartheta-2}}$, $v_{a(t_{\vartheta-1},1)}$, $v_{a(t_{\vartheta-1},2)}$, $\ldots$, $v_{a(t_{\vartheta-1},k-2)}$, $v_{1}\}$, $\mathcal{U}^{'}=\mathcal{U}-e_{t_{\vartheta-1}}+e^{'}_{t_{\vartheta-1}}$.  As Assertion 1, it proved that $\rho(\mathcal{U}^{'})>\rho(\mathcal{U})$, which contradicts the maximality of $\rho(\mathcal{U})$.

{\bf Case 2} $v_{t_{\vartheta}}\neq v_{1}$. Denote by $e_{t_{\vartheta}}=\{v_{t_{\vartheta-1}}$, $v_{a(t_{\vartheta},1)}$, $v_{a(t_{\vartheta},2)}$, $\ldots$, $v_{a(t_{\vartheta},k-2)}$, $v_{t_{\vartheta}}\}$. Let $e^{'}_{t_{\vartheta}}=\{v_{t_{\vartheta-1}}$, $v_{a(t_{\vartheta},1)}$, $v_{a(t_{\vartheta},2)}$, $\ldots$, $v_{a(t_{\vartheta},k-2)}$, $v_{1}\}$, $\mathcal{U}^{'}=\mathcal{U}-e_{t_{\vartheta}}+e^{'}_{t_{\vartheta}}$.  As Assertion 1, it proved that $\rho(\mathcal{U}^{'})>\rho(\mathcal{U})$, which contradicts the maximality of $\rho(\mathcal{U})$.

From the above two cases for Claim 2, Claim 2 is obtained. Moreover, we get that $\mathcal{U}\cong \mathscr{U}^{\ast}$ by above narrations.
 This completes the proof. \ \ \ \ \ $\Box$
\end{Proof}

\

\setlength{\unitlength}{0.6pt}
\begin{center}
\begin{picture}(578,116)
\qbezier(0,60)(0,72)(28,81)\qbezier(28,81)(57,90)(98,90)\qbezier(98,90)(138,90)(167,81)\qbezier(167,81)(196,72)(196,60)\qbezier(196,60)(196,47)(167,38)
\qbezier(167,38)(138,30)(98,30)\qbezier(98,30)(57,30)(28,38)\qbezier(28,38)(0,47)(0,60)
\qbezier(132,60)(132,74)(176,85)\qbezier(176,85)(221,95)(285,95)\qbezier(285,95)(348,95)(393,85)\qbezier(393,85)(438,74)(438,60)\qbezier(438,60)(438,45)(393,34)
\qbezier(393,34)(348,24)(285,24)\qbezier(285,25)(221,25)(176,34)\qbezier(176,34)(131,45)(132,60)
\qbezier(373,60)(373,72)(402,80)\qbezier(402,80)(432,89)(475,89)\qbezier(475,89)(517,89)(547,80)\qbezier(547,80)(577,72)(577,60)\qbezier(577,60)(577,47)(547,39)
\qbezier(547,39)(517,30)(475,30)\qbezier(475,31)(432,31)(402,39)\qbezier(402,39)(372,47)(373,60)
\put(310,65){\circle*{4}}
\put(299,65){\circle*{4}}
\put(164,65){\circle*{4}}
\put(157,51){$v_{1}$}
\put(289,65){\circle*{4}}
\put(252,65){\circle*{4}}
\put(221,65){\circle*{4}}
\put(406,66){\circle*{4}}
\put(399,51){$v_{k}$}
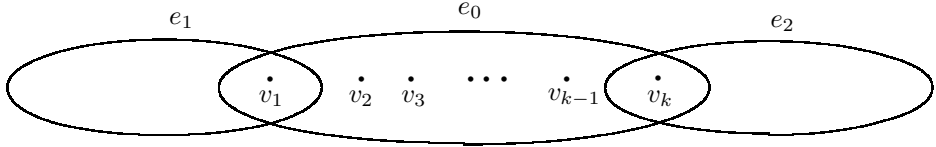
\put(195,-9){Fig. 3.1. $e_{0}$, $e_{1}$, $e_{2}$ in $\mathcal{G}$}
\put(101,100){$e_{1}$}
\put(281,106){$e_{0}$}
\put(476,98){$e_{2}$}
\put(213,51){$v_{2}$}
\put(246,51){$v_{3}$}
\put(349,65){\circle*{4}}
\put(337,52){$v_{k-1}$}
\end{picture}
\end{center}

\begin{lemma}{\bf \cite{GYSLZ}} \label{le3,08}
Let $\mathcal{G}$ be a hypergraph with spectral radius $\rho$, $e_{0}$, $e_{1}$, $e_{2}$ be three edges in $\mathcal{G}$ with $e_{0}=\{v_{1}$, $v_{2}$, $\ldots$, $v_{k-1}$, $v_{k}\}$, satisfying that $deg_{\mathcal{G}}(v_{2})=deg_{\mathcal{G}}(v_{3})=\cdots deg_{\mathcal{G}}(v_{k-1})=1$ ($k\geq 3$), $e_{1}\cap e_{0}=\{v_{1}\}$, $e_{2}\cap e_{0}=\{v_{k}\}$ (see Fig. 3.1.). Let $X$ be
the $principal$ $eigenvector$ of hypergraph $\mathcal{G}$. Then $x_{v_{2}}=x_{v_{3}}=\cdots=x_{v_{k-1}}=\frac{x_{v_{1}}+x_{v_{k}}}{(k-1)\rho-(k-3)}<\min\{x_{v_{1}}, x_{v_{k}}\}$.
\end{lemma}

\begin{lemma}{\bf \cite{GYSLZ}} \label{le3,08,02}
Let $\mathcal{G}$ be a hypergraph with spectral radius $\rho$, $e=\{u$, $v_{1}$, $v_{2}$, $\ldots$, $v_{k-1}\}$ be a pendant edge in $\mathcal{G}$ ($k\geq 2$), where $deg_{\mathcal{G}}(u)\geq 2$. Then in the principal eigenvector $X$ of $\mathcal{G}$, $x_{v_{1}}=x_{v_{2}}=\cdots=x_{v_{k-1}}=\frac{x_{u}}{(k-1)\rho-(k-2)}<x_{u}$.
\end{lemma}

Let $G = (V, E)$ be a $k$-uniform hypergraph, and $e = \{u_{1}$, $u_{2}$, $\ldots$, $u_{k}\}$, $f = \{v_{1}, v_{2}, \ldots, v_{k}\}$
be two edges of $G$, and let $e^{'} = (e \setminus U_{1})\cup V_{1}$, $f^{'} = (f \setminus V_{1}) \cup U_{1}$, $G^{'} = (V, E^{'})$, where $U_{1} = \{u_{1}$, $u_{2}$, $\ldots$, $u_{r}\}$,
$V_{1} = \{v_{1}$, $v_{2}$, $\ldots$, $v_{r}\}$, $1\leq r\leq k-1$, $E^{'} =
(E \setminus\{e, f\})\cup \{e^{'}, f^{'}\}$. Then we denote by $G^{'}=G\langle {u_{1},u_{2},\ldots,u_{r}\rightleftharpoons
{v_{1},v_{2},\ldots,v_{r}}} \rangle$ or $G^{'}=G\langle {U_{1}\rightleftharpoons
{V_{1}}} \rangle$.
Let $X$ be an eigenvector on a connected $k$-uniform hypergraph $G$. For
the simplicity, we let $x_{S} =\sum_{v\in S} x_{v}$ where $S\subseteq V$, and for an edge $e$, we let $x_{e} =\sum_{v\in e} x_{v}$.

\begin{lemma} \label{le3,42}
Let $G = (V, E)$ be a connected $k$-uniform hypergraph, and $e = \{u_{1}$, $u_{2}$, $\ldots$, $u_{k}\}$, $f = \{v_{1}$, $v_{2}$, $\ldots$, $v_{k}\}$
be two edges of $G$, $U_{1} = \{u_{1}$, $u_{2}$, $\ldots$, $u_{r}\}$,
$V_{1} = \{v_{1}$, $v_{2}$, $\ldots$, $v_{r}\}$, $1\leq r\leq k-1$. $G^{'}=G\langle{U_{1}\rightleftharpoons
{V_{1}}}  \rangle$. Let $X$ be the principal eigenvector on $G$, $U_{2} = e \setminus U_{1}$, $V_{2} = f \setminus V_{1}$.  If $x_{U_{1}} \geq x_{V_{1}}$, $x_{U_{2}} \leq x_{V_{2}}$, and $G^{'}$ is connected, then
$\rho(G) \leq \rho(G^{'})$. Moreover, we have

(1) if $x_{U_{1}} > x_{V_{1}}$, $x_{U_{2}} < x_{V_{2}}$, then $\rho(G) < \rho(G^{'})$;

(2) if $\rho(G) = \rho(G^{'})$, then $X$ is also the principal eigenvector of $G^{'}$.
\end{lemma}

\begin{proof}
Note that
$$\rho(G^{'})- \rho(G)\geq X^{T}\mathcal {A}(G^{'})X-X^{T}\mathcal {A}(G)X=\frac{2}{k-1}(x_{U_{1}}-x_{V_{1}})(x_{V_{2}}-x_{U_{2}}).\hspace{2cm} (\ast)$$
Thus it follows that if $x_{U_{1}} \geq x_{V_{1}}$, $x_{U_{2}} \leq x_{V_{2}}$, then $\rho(G^{'})\geq \rho(G)$. Moreover, if $x_{U_{1}} > x_{V_{1}}$, $x_{U_{2}} < x_{V_{2}}$, then $\rho(G) < \rho(G^{'})$; if $\rho(G) = \rho(G^{'})$, by the preliminary in Section 2, then it follows that $X$ is also the principal eigenvector of $G^{'}$.
This completes the proof. \ \ \ \ \ $\Box$
\end{proof}

\setlength{\unitlength}{0.7pt}
\begin{center}
\begin{picture}(535,486)
\put(369,350){\circle*{4}}
\put(504,350){\circle*{4}}
\qbezier(369,350)(440,486)(504,350)
\qbezier(369,350)(438,404)(504,350)
\qbezier(369,350)(435,239)(504,350)
\qbezier(369,350)(439,306)(504,350)
\put(34,144){\circle*{4}}
\qbezier(369,350)(297,374)(294,349)
\qbezier(294,349)(296,321)(369,350)
\put(36,343){\circle*{4}}
\put(179,343){\circle*{4}}
\qbezier(36,343)(114,467)(179,343)
\qbezier(36,343)(111,385)(179,343)
\qbezier(36,343)(114,305)(179,343)
\qbezier(36,343)(115,235)(179,343)
\put(106,384){\circle*{4}}
\qbezier(26,400)(14,370)(106,384)
\qbezier(26,400)(32,424)(106,384)
\qbezier(106,384)(36,434)(59,449)
\qbezier(59,449)(84,465)(106,384)
\qbezier(106,384)(106,467)(131,458)
\qbezier(106,384)(156,453)(131,458)
\qbezier(106,384)(188,421)(186,395)
\qbezier(106,384)(191,371)(186,395)
\put(158,428){\circle*{4}}
\put(153,436){\circle*{4}}
\put(163,419){\circle*{4}}
\put(30,134){\circle*{4}}
\put(97,259){$\mathcal {F}_{1}$}
\put(428,264){$\mathcal {F}_{2}$}
\put(171,-9){Fig. 3.2. $\mathcal {F}_{1}-\mathcal {F}_{3}$, $\mathcal {F}_{(\mathcal {R}; \mathcal {S}; \mathcal {T})}$}
\put(150,362){$e_{2}$}
\put(183,341){$v_{2}$}
\put(135,312){$e_{1}$}
\put(30,394){$e_{3}$}
\put(61,436){$e_{4}$}
\put(22,342){$v_{1}$}
\put(508,348){$v_{2}$}
\put(306,348){$e_{3}$}
\put(26,124){\circle*{4}}
\put(359,357){$v_{1}$}
\put(458,315){$e_{1}$}
\put(463,380){$e_{2}$}
\put(99,372){$v_{\eta}$}
\put(367,411){$e_{4}$}
\put(437,394){\circle*{4}}
\qbezier(358,414)(365,442)(437,394)
\qbezier(358,414)(351,389)(437,394)
\qbezier(437,394)(373,440)(396,456)
\qbezier(396,456)(422,467)(437,394)
\qbezier(437,394)(509,394)(508,416)
\qbezier(437,394)(506,444)(508,416)
\put(441,441){\circle*{4}}
\put(452,439){\circle*{4}}
\put(463,436){\circle*{4}}
\put(119,444){$e_{5}$}
\put(397,444){$e_{5}$}
\put(162,391){$e_{m}$}
\put(485,411){$e_{m}$}
\put(68,115){\circle*{4}}
\put(197,114){\circle*{4}}
\qbezier(68,115)(136,177)(197,114)
\qbezier(68,115)(136,255)(197,114)
\qbezier(68,115)(137,76)(197,114)
\qbezier(68,115)(137,7)(197,114)
\qbezier(68,115)(80,52)(53,54)
\qbezier(68,115)(36,69)(53,54)
\qbezier(68,115)(25,73)(13,91)
\qbezier(68,115)(0,119)(13,91)
\qbezier(68,115)(74,172)(51,170)
\qbezier(51,170)(25,163)(68,115)
\put(133,165){\circle*{4}}
\qbezier(133,165)(108,227)(133,226)
\qbezier(133,226)(159,228)(133,165)
\put(76,113){$v_{1}$}
\put(202,112){$v_{2}$}
\put(152,82){$e_{1}$}
\put(148,160){$e_{2}$}
\put(54,71){$e_{3}$}
\put(19,98){$e_{4}$}
\put(36,175){$e_{m-1}$}
\put(125,210){$e_{m}$}
\put(124,29){$\mathcal {F}_{3}$}
\put(431,383){$v_{\eta}$}
\put(127,156){$v_{\eta}$}
\put(347,135){\circle*{4}}
\put(478,135){\circle*{4}}
\qbezier(347,135)(414,273)(478,135)
\qbezier(347,135)(414,193)(478,135)
\qbezier(347,135)(415,25)(478,135)
\qbezier(347,135)(415,92)(478,135)
\qbezier(347,135)(329,78)(305,89)
\qbezier(347,135)(283,113)(305,89)
\qbezier(347,135)(327,195)(304,175)
\qbezier(304,175)(287,155)(347,135)
\put(367,64){\circle*{4}}
\put(356,70){\circle*{4}}
\put(362,67){\circle*{4}}
\put(359,211){\circle*{4}}
\put(349,201){\circle*{4}}
\put(354,206){\circle*{4}}
\put(303,124){\circle*{4}}
\put(303,141){\circle*{4}}
\put(303,133){\circle*{4}}
\qbezier(478,135)(499,190)(519,174)
\qbezier(478,135)(535,154)(519,174)
\qbezier(478,135)(487,79)(513,88)
\qbezier(478,135)(531,112)(513,88)
\put(515,122){\circle*{4}}
\put(515,139){\circle*{4}}
\put(515,131){\circle*{4}}
\put(382,169){\circle*{4}}
\qbezier(382,169)(321,167)(334,191)
\qbezier(334,191)(342,205)(382,169)
\qbezier(382,169)(395,225)(374,221)
\qbezier(374,221)(354,221)(382,169)
\put(381,106){\circle*{4}}
\qbezier(381,106)(331,96)(336,78)
\qbezier(381,106)(349,58)(336,78)
\qbezier(381,106)(404,65)(385,56)
\qbezier(381,106)(363,52)(385,56)
\put(443,213){\circle*{4}}
\put(420,213){\circle*{4}}
\put(432,213){\circle*{4}}
\put(442,65){\circle*{4}}
\put(418,65){\circle*{4}}
\put(430,65){\circle*{4}}
\put(432,101){$e_{1}$}
\put(435,173){$e_{2}$}
\put(354,132){$v_{1}$}
\put(458,134){$v_{2}$}
\put(381,31){$\mathcal {F}_{(\mathcal {R}; \mathcal {S}; \mathcal {T})}$}
\end{picture}
\end{center}

Let $\mathcal {F}_{1}$ be the $k$-uniform unicyclic hypergraph consisting of hypercycle $\mathbb{C}$ (narrated in Section 1) and $\frac{n}{k-1}-2$ pendant edges attaching at vertex $\eta\in (e_{2}\setminus\{v_{1}, v_{2}\})$ (see Fig. 3.2); let $\mathcal {F}_{2}$ be the $k$-uniform unicyclic hypergraph consisting of hypercycle $\mathbb{C}$, $\frac{n}{k-1}-3$ pendant edges attaching at vertex $\eta\in (e_{2}\setminus\{v_{1}, v_{2}\})$ and one pendant edge attaching at vertex $v_{1}$ (see Fig. 3.2); let $\mathcal {F}_{3}$ be the $k$-uniform unicyclic hypergraph consisting of hypercycle $\mathbb{C}$, $\frac{n}{k-1}-3$ pendant edges attaching at vertex $v_{1}$ and one pendant edge attaching at vertex $\eta\in (e_{2}\setminus\{v_{1}, v_{2}\})$ (see Fig. 3.2). Let $\mathcal {R}=[r_{1}, r_{2}]$, $\mathcal {S}=[s_{1}, s_{2}, \ldots, s_{k-2}]$, $\mathcal {T}=[t_{1}, t_{2}, \ldots, t_{k-2}]$ be three multisets (or arrays) satisfying that $r_{1}+ r_{2}+\sum^{k-2}_{i=1}s_{i}+\sum^{k-2}_{i=1}t_{i}=\frac{n}{k-1}-2$ where $r_{i}$, $s_{i}$, $t_{i}$ may be zero for some $i$, and it is possible that $r_{1}=r_{2}$, $s_{i}=s_{j}$ for $i\neq j$, $t_{i}=t_{j}$ for $i\neq j$; $\mathcal {F}_{(\mathcal {R}; \mathcal {S}; \mathcal {T})}\in \Lambda$ be the $k$-uniform unicyclic hypergraph which consists of hypercycle $\mathbb{C}$, $r_{i}$ pendant edges attaching at $v_{i}$ for $i=1, 2$, $s_{i}$ pendant edges attaching at $v_{a(1, i)}$ for $i=1, 2, \ldots, k-2$, $t_{i}$ pendant edges attaching at $v_{a(2, i)}$ for $i=1, 2, \ldots, k-2$ (see Fig. 3.2).

\begin{lemma}\label{le01.01.01} 
$\rho(\mathcal {F}_{1})<\rho(\mathcal {F})$.
\end{lemma}

\begin{proof}
In $\mathcal {F}_{1}$, denote by $e_{1}=\{v_{1}$, $v_{a(1,1)}$, $v_{a(1,2)}$, $\ldots$, $v_{a(1,k-2)}$, $v_{2}\}$, $e_{2}=\{v_{1}$, $v_{a(2,1)}$, $v_{a(2,2)}$, $\ldots$, $v_{a(2,k-3)}$, $v_{\eta}$, $v_{2}\}$, $e_{i}=\{v_{\eta}$, $v_{a(i,1)}$, $v_{a(i,2)}$, $\ldots$, $v_{a(i,k-1)}\}$ for $3\leq i\leq m$, where $m=\frac{n}{k-1}$. For short, we let $\rho_{1}=\rho(\mathcal {F}_{1})$. Let $X$ be
the principal eigenvector of $\mathcal {F}_{1}$. By Lemma \ref{le03,01}, we get that $x_{v_{1}}=x_{v_{2}}$, $x_{v_{a(2,1)}}=x_{v_{a(2,2)}}=\cdots =v_{v_{a(2,k-3)}}$. By Lemma \ref{le3,08} and Lemma \ref{le3,08,02}, we have $x_{v_{a(1,i)}}=\frac{2x_{v_{1}}}{(k-1)\rho_{1}-(k-3)}<x_{v_{1}}$ for $1\leq i\leq k-2$, $x_{v_{a(j,i)}}=\frac{x_{v_{\eta}}}{(k-1)\rho_{1}-(k-2)}<x_{v_{\eta}}$ for $3\leq j\leq m$, $1\leq i\leq k-1$.

{\bf Assertion 1} $x_{v_{a(2,1)}}<x_{v_{\eta}}$. Otherwise, suppose $x_{v_{a(2,1)}}\geq x_{v_{\eta}}$. Let $e^{'}_{i}=\{v_{a(2,1)}$, $v_{a(i,1)}$, $v_{a(i,2)}$, $\ldots$, $v_{a(i,k-1)}\}$ for $3\leq i\leq m$, and $\mathcal {F}^{'}_{1}=\mathcal {F}_{1}-\sum_{3\leq i\leq m}e_{i}+\sum_{3\leq i\leq m}e^{'}_{i}$. Using Lemma \ref{le03.02} gets that $\rho(\mathcal {F}^{'}_{1})>\rho(\mathcal {F}_{1})$, which causes a contradiction because $\mathcal {F}^{'}_{1}\cong \mathcal {F}_{1}$.

{\bf Assertion 2} $x_{v_{a(2,1)}}<x_{v_{1}}$. Otherwise, suppose $x_{v_{a(2,1)}}\geq x_{v_{1}}$. Let $e^{'}_{1}=\{v_{2}$, $v_{a(1,1)}$, $v_{a(1,2)}$, $\ldots$, $v_{a(1,k-2)}$, $v_{a(2,1)}\}$, and $\mathcal {F}^{'}_{1}=\mathcal {F}_{1}-e_{1}+e^{'}_{1}$. As Assertion 1, we gets that $\rho(\mathcal {F}^{'}_{1})>\rho(\mathcal {F}_{1})$, which causes a contradiction because $\mathcal {F}^{'}_{1}\cong \mathcal {F}_{1}$.

{\bf Assertion 3} $x_{v_{a(2,i)}}=\frac{2x_{v_{1}}+x_{v_{\eta}}}{(k-1)\rho_{1}-(k-4)}<\min\{x_{v_{1}}, x_{v_{\eta}}\}$ for $1\leq i\leq k-3$. This assertion follows from $x_{v_{1}}=x_{v_{2}}$, $x_{v_{a(2,i)}}=x_{v_{a(2,j)}}$ for $1\leq i< j\leq k-3$, Assertion 1 and Assertion 2, and $\rho_{1} x_{v_{a(2,1)}}=\frac{1}{k-1}((k-4)x_{v_{a(2,1)}}+2x_{v_{1}}+x_{v_{\eta}})$.

To prove our lemma, we need consider two cases further.

{\bf Case 1} $x_{v_{1}}\geq x_{v_{\eta}}$. Let $e^{'}_{3}=\{v_{2}$, $v_{a(3,1)}$, $v_{a(3,2)}$, $\ldots$, $v_{a(3,k-1)}\}$, $e^{'}_{i}=\{v_{1}$, $v_{a(i,1)}$, $v_{a(i,2)}$, $\ldots$, $v_{a(i,k-1)}\}$ for $4\leq i\leq m$, and $\mathcal {F}^{'}_{1}=\mathcal {F}_{1}-\sum_{3\leq i\leq m}e_{i}+\sum_{3\leq i\leq m}e^{'}_{i}$. Note that $x_{v_{1}}=x_{v_{2}}$. Then $$X^{T}\mathcal {A}_{\mathcal {F}^{'}_{1}}X-X^{T}\mathcal {A}_{\mathcal {F}_{1}}X=\frac{2}{k-1}(\sum^{k-1}_{i=1}x_{v_{a(3,i)}})(x_{v_{2}}-x_{v_{\eta}})+\frac{2}{k-1}(\sum^{m}_{j=4}\sum^{k-1}_{i=1}x_{v_{a(j,i)}})(x_{v_{1}}-x_{v_{\eta}})$$
$$ =\frac{2}{k-1}(\sum^{m}_{j=3}\sum^{k-1}_{i=1}x_{v_{a(j,i)}})(x_{v_{1}}-x_{v_{\eta}})\geq 0.$$ Thus, we get that $\rho(\mathcal {F}^{'}_{1})\geq\rho(\mathcal {F}_{1})$.
If $\rho(\mathcal {F}^{'}_{1})=\rho(\mathcal {F}_{1})$, then $\rho(\mathcal {F}^{'}_{1})=X^{T}\mathcal {A}_{\mathcal {F}^{'}_{1}}X=X^{T}\mathcal {A}_{\mathcal {F}_{1}}X=\rho(\mathcal {F}_{1}).$ Then $X$ is also an eigenvector corresponding to $\rho(\mathcal {F}^{'}_{1})$. That is $\mathcal {A}_{\mathcal {F}^{'}_{1}}X=\rho(\mathcal {F}^{'}_{1})X$. But $(\mathcal {A}_{\mathcal {F}^{'}_{1}}X)_{v_{2}}-(\mathcal {A}_{\mathcal {F}_{1}}X)_{v_{2}}=\frac{1}{k-1}\sum^{k-1}_{i=1}x_{v_{a(3,i)}}> 0$ which contradicts $\mathcal {A}_{\mathcal {F}^{'}_{1}}X=\rho(\mathcal {F}^{'}_{1})X$. Then it follows that $\mathcal {A}_{\mathcal {F}^{'}_{1}}X>\rho(\mathcal {F}^{'}_{1})X$. Note that $\mathcal {F}^{'}_{1}\cong \mathcal {F}$ now. Thus the lemma follows as desired.

{\bf Case 2} $x_{v_{1}}< x_{v_{\eta}}$.

{\bf Claim 1} $x_{v_{a(2,i)}}\geq x_{v_{a(j,s)}}$ for $1\leq i\leq k-3$, $1\leq s\leq k-1$, $3\leq j\leq m$. We prove this claim by contradiction. Otherwise, suppose $x_{v_{a(2,i)}}<x_{v_{a(3,1)}}$.

{\bf Subcase 2.1} $x_{v_{1}}> x_{v_{a(3,1)}}$. Let $\mathbb{A}_{1}=\{v_{a(3,1)}$, $v_{a(3,2)}$, $\ldots$, $v_{a(3,k-3)}\}$, $\mathbb{B}_{1}=\{v_{a(2,1)}$, $v_{a(2,2)}$, $\ldots$, $v_{a(2,k-3)}\}$, $\mathcal {F}^{'}_{1}=F_{1}\langle
{\mathbb{A}_{1}\rightleftharpoons
{\mathbb{B}_{1}}} \rangle$. By Lemma \ref{le3,42}, we get that $\rho(\mathcal {F}^{'}_{1})>\rho(\mathcal {F}_{1})$, which causes a contradiction because $\mathcal {F}^{'}_{1}\cong \mathcal {F}_{1}$.

{\bf Subcase 2.2} $x_{v_{1}}\leq x_{v_{a(3,1)}}$. Let $\mathbb{A}_{1}=\{v_{a(3,1)}$, $v_{a(3,2)}$, $\ldots$, $v_{a(3,k-1)}\}$, $\mathbb{B}_{1}=\{v_{a(2,1)}$, $v_{a(2,2)}$, $\ldots$, $v_{a(2,k-3)}$, $v_{1}$, $v_{2}\}$, $\mathcal {F}^{'}_{1}=F_{1}\langle
{\mathbb{A}_{1}\rightleftharpoons
{\mathbb{B}_{1}}} \rangle$. By Lemma \ref{le3,42}, we get that $\rho(\mathcal {F}^{'}_{1})\geq\rho(\mathcal {F}_{1})$. If $\rho(\mathcal {F}^{'}_{1})=\rho(\mathcal {F}_{1})$, using Lemma \ref{le3,42} again gets that $\mathcal {A}_{\mathcal {F}^{'}_{1}}X=\rho(\mathcal {F}_{1})X$. Suppose $\rho(\mathcal {F}^{'}_{1})=\rho(\mathcal {F}_{1})$. Then $(\mathcal {A}_{\mathcal {F}^{'}_{1}}X)_{v_{a(3,i)}}=(\mathcal {A}_{\mathcal {F}_{1}}X)_{v_{a(3,i)}}$ for $1\leq i\leq k-1$. Suppose in $\mathcal {F}^{'}_{1}$, $e^{'}_{1}=\{v_{a(3,1)}$, $v_{a(1,1)}$, $v_{a(1,2)}$, $\ldots$, $v_{a(1,k-2)}$, $v_{a(3,2)}\}$, $e^{'}_{2}=\{v_{a(3,1)}$, $v_{a(2,1)}$, $v_{a(2,2)}$, $\ldots$, $v_{a(2,k-3)}$, $v_{\eta}$, $v_{a(3,2)}\}$. But $(\mathcal {A}_{\mathcal {F}^{'}_{1}}X)_{v_{a(3,1)}}>(\mathcal {A}_{\mathcal {F}_{1}}X)_{v_{a(3,1)}}$ which contradicts $(\mathcal {A}_{\mathcal {F}^{'}_{1}}X)_{v_{a(3,1)}}=(\mathcal {A}_{\mathcal {F}_{1}}X)_{v_{a(3,1)}}$.

Above two subcases means that $x_{v_{a(2,i)}}\geq x_{v_{a(j,s)}}$ for $1\leq i\leq k-3$, $1\leq s\leq k-1$, $3\leq j\leq m$. Then Claim 1 follows.

Combing Assertion 2 and Claim 1, we get the following Claim 2.

{\bf Claim 2} $x_{v_{1}}> x_{v_{a(j,s)}}$ for $1\leq s\leq k-1$, $3\leq j\leq m$.

Thus $x_{v_{1}}> x_{v_{a(3,1)}}$. Note that $x_{v_{1}}< x_{v_{\eta}}$.

{\bf Subcase $2.1^{'}$} $x_{v_{a(1,1)}}\geq x_{v_{a(3,1)}}$. Let $e^{'}_{1}=\{v_{2}$, $v_{a(1,1)}$, $v_{a(1,2)}$, $\ldots$, $v_{a(1,k-2)}$, $v_{\eta}\}$, $e^{'}_{3}=\{v_{2}$, $v_{a(3,1)}$, $v_{a(3,2)}$, $\ldots$, $v_{a(3,k-2)}$, $v_{a(3,k-1)}\}$, $\mathcal {F}^{'}_{1}=\mathcal {F}_{1}-e_{1}-e_{3}+e^{'}_{1}+e^{'}_{3}$.  Note that $x_{v_{1}}=x_{v_{2}}$. Then
$$X^{T}\mathcal {A}_{\mathcal {F}^{'}_{1}}X-X^{T}\mathcal {A}_{\mathcal {F}_{1}}X=\frac{2}{k-1}(x_{v_{2}}+\sum^{k-2}_{i=1}x_{v_{a(1,i)}})(x_{v_{\eta}}-x_{v_{1}})+\frac{2}{k-1}(\sum^{k-1}_{i=1}x_{v_{a(3,i)}})(x_{v_{2}}-x_{v_{\eta}})$$$$=\frac{2}{k-1}(x_{v_{2}}+\sum^{k-2}_{i=1}x_{v_{a(1,i)}}-\sum^{k-1}_{i=1}x_{v_{a(3,i)}})(x_{v_{\eta}}-x_{v_{1}})> 0.$$ Note that $\mathcal {F}^{'}_{1}\cong \mathcal {F}$.
It follows that $\rho(\mathcal {F})=\rho(\mathcal {F}^{'}_{1})>\rho(\mathcal {F}_{1})$.

{\bf Subcase $2.2^{'}$} $x_{v_{a(1,1)}}< x_{v_{a(3,1)}}$.
Let $e^{'}_{1}=\{v_{2}$, $v_{a(1,1)}$, $v_{a(1,2)}$, $\ldots$, $v_{a(1,k-2)}$, $v_{a(3,k-1)}\}$, $e^{'}_{3}=\{v_{2}$, $v_{a(3,1)}$, $v_{a(3,2)}$, $\ldots$, $v_{a(3,k-2)}$, $v_{\eta}\}$, $\mathcal {F}^{'}_{1}=\mathcal {F}_{1}-e_{1}-e_{3}+e^{'}_{1}+e^{'}_{3}$. Note that $x_{v_{1}}=x_{v_{2}}$. Then
$$X^{T}\mathcal {A}_{\mathcal {F}^{'}_{1}}X-X^{T}\mathcal {A}_{\mathcal {F}_{1}}X=\frac{2}{k-1}(x_{v_{\eta}}+
\sum^{k-1}_{i=1}x_{v_{a(3,i)}})(x_{v_{2}}-x_{v_{a(3,k-1)}})+\frac{2}{k-1}(x_{v_{2}}+\sum^{k-2}_{i=1}x_{v_{a(1,i)}})(x_{v_{a(3,k-1)}}-x_{v_{1}})$$
$$=\frac{2}{k-1}(x_{v_{\eta}}+
\sum^{k-2}_{i=1}x_{v_{a(3,i)}}-x_{v_{2}}-\sum^{k-2}_{i=1}x_{v_{a(1,i)}})(x_{v_{2}}-x_{v_{a(3,k-1)}})> 0.$$
Note that  $\mathcal {F}^{'}_{1}\cong \mathcal {F}$. Then it follows that $\rho(\mathcal {F})=\rho(\mathcal {F}^{'}_{1})>\rho(\mathcal {F}_{1})$.

From the above narrations, we conclude that $\rho(\mathcal {F})>\rho(\mathcal {F}_{1})$. This completes the proof. \ \ \ \ \ $\Box$
\end{proof}

\begin{lemma}\label{le01.01.01,00} 
$\rho(\mathcal {F}_{i})<\rho(\mathcal {F})$ for $i=2, 3$.
\end{lemma}

\begin{proof}
Let $X$ be
the principal eigenvector of $\mathcal {F}_{2}$. In $\mathcal {F}_{2}$, denote by $e_{1}=\{v_{1}$, $v_{a(1,1)}$, $v_{a(1,2)}$, $\ldots$, $v_{a(1,k-2)}$, $v_{2}\}$, $e_{2}=\{v_{1}$, $v_{a(2,1)}$, $v_{a(2,2)}$, $\ldots$, $v_{a(2,k-3)}$, $v_{\eta}$, $v_{2}\}$, $e_{3}=\{v_{1}$, $v_{a(3,1)}$, $v_{a(3,2)}$, $\ldots$, $v_{a(3,k-1)}\}$, $e_{i}=\{v_{\eta}$, $v_{a(i,1)}$, $v_{a(i,2)}$, $\ldots$, $v_{a(i,k-1)}\}$ for $4\leq i\leq m$, where $m=\frac{n}{k-1}$. If $x_{v_{\eta}}\geq v_{2}$, then we let $e^{'}_{1}=\{v_{1}$, $v_{a(1,1)}$, $v_{a(1,2)}$, $\ldots$, $v_{a(1,k-2)}$, $v_{\eta}\}$, and $\mathcal {F}^{'}_{2}=\mathcal {F}_{2}-e_{1}+e^{'}_{1}$; if $x_{v_{\eta}}< v_{2}$, then we let $e^{'}_{i}=\{v_{2}$, $v_{a(i,1)}$, $v_{a(i,2)}$, $\ldots$, $v_{a(i,k-1)}\}$ for $4\leq i\leq m$, and $\mathcal {F}^{'}_{2}=\mathcal {F}_{2}-\sum^{m}_{i=4}e_{i}+\sum^{m}_{i=4}e^{'}_{i}$. Using Lemma \ref{le03.02} gets that $\rho(\mathcal {F}_{2})<\rho(\mathcal {F}^{'}_{2})$. Note that $\mathcal {F}^{'}_{2}\cong \mathcal {F}$. Then we have $\rho(\mathcal {F}_{2})<\rho(\mathcal {F})$.

In a same way, we get that $\rho(\mathcal {F}_{3})<\rho(\mathcal {F})$.
This completes the proof. \ \ \ \ \ $\Box$
\end{proof}

\begin{lemma}\label{le01.01.01,01} 
For $\mathcal {F}_{(\mathcal {R}; \mathcal {S}; \mathcal {T})}\in \Lambda$, we have $\rho(\mathcal {F}_{(\mathcal {R}; \mathcal {S}; \mathcal {T})})<\rho(\mathcal {F})$.
\end{lemma}

\begin{proof}
In $F_{(\mathcal {R}; \mathcal {S}; \mathcal {T})}$, denote by $e_{1}=\{v_{1}$, $v_{a(1,1)}$, $v_{a(1,2)}$, $\ldots$, $v_{a(1,k-2)}$, $v_{2}\}$, $e_{2}=\{v_{1}$, $v_{a(2,1)}$, $v_{a(2,2)}$, $\ldots$, $v_{a(2,k-3)}$, $v_{a(2,k-2)}$, $v_{2}\}$; $e_{(\kappa _{i},j)}=\{v_{i}$, $v_{a(\kappa _{i},j,1)}$, $v_{a(\kappa _{i},j,2)}$, $\ldots$, $v_{a(\kappa _{i},j,k-1)}\}$ the $jth$ pendant edge attaching at vertex $v_{i}$ where $1\leq i\leq 2$, $r_{i}\geq 1$, $j=1, 2, \ldots, r_{i}$; $e_{(\varepsilon_{i},j)}=\{v_{a(1,i)}$, $v_{a(\varepsilon_{i},j,1)}$, $v_{a(\varepsilon _{i},j,2)}$, $\ldots$, $v_{a(\varepsilon _{i},j,k-1)}\}$ the $jth$ pendant edge attaching at vertex $v_{a(1,i)}$ where $1\leq i\leq k-2$, $s_{i}\geq 1$, $j=1, 2, \ldots, s_{i}$; $e_{(\vartheta _{i},j)}=\{v_{a(2,i)}$, $v_{a(\vartheta _{i},j,1)}$, $v_{a(\vartheta _{i},j,2)}$, $\ldots$, $v_{a(\vartheta_{i},j,k-1)}\}$ the $jth$ pendant edge attaching at vertex $v_{a(2,i)}$ where $1\leq i\leq k-2$, $t_{i}\geq 1$, $j=1, 2, \ldots, t_{i}$.

Let $X$ be
the principal eigenvector of $\mathcal {F}_{(\mathcal {R}; \mathcal {S}; \mathcal {T})}$, $D=\{x_{v_{a(1,i)}}|\, s_{i}\geq 1\}\cup \{x_{v_{a(2,i)}}|\, t_{i}\geq 1\}$. Without loss of generality, we suppose $s_{1}\geq 1$, and $x_{v_{a(1,1)}}=\max\{x_{v}|\, v\in D\}$. Let $e'_{(\varepsilon_{i},j)}=\{v_{a(1,1)}$, $v_{a(\varepsilon_{i},j,1)}$, $v_{a(\varepsilon _{i},j,2)}$, $\ldots$, $v_{a(\varepsilon _{i},j,k-1)}\}$ where $2\leq i\leq k-2$, $s_{i}\geq 1$, $j=1, 2, \ldots, s_{i}$;  $e'_{(\vartheta _{i},j)}=\{v_{a(1,1)}$, $v_{a(\vartheta _{i},j,1)}$, $v_{a(\vartheta _{i},j,2)}$, $\ldots$, $v_{a(\vartheta_{i},j,k-1)}\}$ where $1\leq i\leq k-2$, $t_{i}\geq 1$, $j=1, 2, \ldots, t_{i}$. Let
$$\mathcal {F}'_{(\mathcal {R}; \mathcal {S}; \mathcal {T})}=\mathcal {F}_{(\mathcal {R}; \mathcal {S}; \mathcal {T})}-\sum_{2\leq i\leq k-2, s_{i}\geq 1}\sum^{s_{i}}_{j=1}e_{(\varepsilon_{i},j)}+\sum_{2\leq i\leq k-2, s_{i}\geq 1}\sum^{s_{i}}_{j=1}e'_{(\varepsilon_{i},j)}\hspace{0.2cm}$$$$\hspace{0.4cm}-\sum_{1\leq i\leq k-2, t_{i}\geq 1}\sum^{t_{i}}_{j=1}e_{(\vartheta_{i},j)}+\sum_{1\leq i\leq k-2, t_{i}\geq 1}\sum^{t_{i}}_{j=1}e'_{(\vartheta_{i},j)}.$$
As Case 1 in the proof of Lemma \ref{le01.01.01}, we get that $\rho(\mathcal {F}'_{(\mathcal {R}; \mathcal {S}; \mathcal {T})})>\mathcal {F}_{(\mathcal {R}; \mathcal {S}; \mathcal {T})}$.

Thus hereafter, if there exists some $s_{i}\geq 1$ or some $t_{i}\geq 1$ where $1\leq i\leq k-2$, then we suppose $s_{1}\geq 1$, $s_{i}=0$ for all $2\leq i\leq k-2$, and $t_{i}=0$ for all $1\leq i\leq k-2$.

{\bf Case 1} $r_{1}\geq 1$, $r_{2}\geq 1$.

{\bf Subcase 1.1} In $X$, $\max\{x_{v_{1}}$, $x_{v_{2}}\}\geq x_{v_{a(1,1)}}$. Suppose $\max\{x_{v_{1}}$, $x_{v_{2}}\}=x_{v_{1}}$. Now, let $e^{'}_{(\varepsilon_{1},j)}=\{v_{1}$, $v_{a(\varepsilon_{1},j,1)}$, $v_{a(\varepsilon _{1},j,2)}$, $\ldots$, $v_{a(\varepsilon _{1},j,k-1)}\}$ where $j=1, 2, \ldots, s_{1}$. If $r_{2}\geq 2$, then we let $e^{'}_{(\kappa _{2},j)}=\{v_{1}$, $v_{a(\kappa _{2},j,1)}$, $v_{a(\kappa _{2},j,2)}$, $\ldots$, $v_{a(\kappa _{2},j,k-1)}\}$ where $j=2, \ldots,  r_{2}$, $$B=\mathcal {F}_{(\mathcal {R}; \mathcal {S}; \mathcal {T})}-\sum^{r_{2}}_{j=2}e_{(\kappa _{2},j)}+\sum^{r_{2}}_{j=2}e^{'}_{(\kappa _{2},j)}-\sum^{s_{1}}_{j=1}e_{(\varepsilon_{1},j)}+\sum^{s_{1}}_{j=1}e^{'}_{(\varepsilon_{1},j)}.$$ Let $$B=\mathcal {F}_{(\mathcal {R}; \mathcal {S}; \mathcal {T})}-\sum^{s_{1}}_{j=1}e_{(\varepsilon_{1},j)}+\sum^{s_{1}}_{j=1}e^{'}_{(\varepsilon_{1},j)}$$ if $r_{2}=1$. As Case 1 in the proof of Lemma \ref{le01.01.01}, we get that $\rho(B)>\rho(\mathcal {F}_{(\mathcal {R}; \mathcal {S}; \mathcal {T})})$. Note that $B\cong \mathcal {F}$. Then $\rho(\mathcal {F})>\rho(\mathcal {F}_{(\mathcal {R}; \mathcal {S}; \mathcal {T})})$.

{\bf Subcase 1.2} In $X$, $\max\{x_{v_{1}}$, $x_{v_{2}}\}< x_{v_{a(1,1)}}$. Suppose $\max\{x_{v_{1}}$, $x_{v_{2}}\}=x_{v_{1}}$. Now, let $e^{'}_{(\kappa _{2},j)}=\{v_{a(1,1)}$, $v_{a(\kappa _{2},j,1)}$, $v_{a(\kappa _{2},j,2)}$, $\ldots$, $v_{a(\kappa _{1},j,k-1)}\}$ where $j=1, \ldots,  r_{2}$; $e^{'}_{2}=\{v_{1}$, $v_{a(2,1)}$, $v_{a(2,2)}$, $\ldots$, $v_{a(2,k-3)}$, $v_{a(2,k-2)}$, $v_{v_{a(1,1)}}\}$. If $r_{1}\geq 2$, then we let $e^{'}_{(\kappa _{1},j)}=\{v_{a(1,1)}$, $v_{a(\kappa _{1},j,1)}$, $v_{a(\kappa _{1},j,2)}$, $\ldots$, $v_{a(\kappa _{2},j,k-1)}\}$ where $j=2, \ldots,  r_{1}$, $$B=\mathcal {F}_{(\mathcal {R}; \mathcal {S}; \mathcal {T})}-\sum^{r_{1}}_{j=2}e_{(\kappa _{1},j)}+\sum^{r_{1}}_{j=2}e^{'}_{(\kappa _{1},j)}-\sum^{r_{2}}_{j=1}e_{(\kappa _{2},j)}+\sum^{r_{2}}_{j=1}e^{'}_{(\kappa _{2},j)}-e_{2}+e^{'}_{2}.$$ Let $$B=\mathcal {F}_{(\mathcal {R}; \mathcal {S}; \mathcal {T})}-\sum^{r_{2}}_{j=1}e_{(\kappa _{2},j)}+\sum^{r_{2}}_{j=1}e^{'}_{(\kappa _{2},j)}-e_{2}+e^{'}_{2}$$ if $r_{1}=1$. As Subcase 1.1, we get that $\rho(B)>\rho(\mathcal {F}_{(\mathcal {R}; \mathcal {S}; \mathcal {T})})$. Note that $B\cong \mathcal {F}$. Then $\rho(\mathcal {F})>\rho(\mathcal {F}_{(\mathcal {R}; \mathcal {S}; \mathcal {T})})$.

{\bf Case 2} Only one of $r_{1}$, $r_{2}$ is zero. Note that $\mathcal {F}_{(\mathcal {R}; \mathcal {S}; \mathcal {T})}\in \Lambda$. Then $s_{1}>0$. Without loss of generality, suppose $r_{2}=0$.

{\bf Subcase 1.1} In $X$, $x_{v_{1}}< x_{v_{a(1,1)}}$. Now, if $r_{1}\geq 2$, then we let $e^{'}_{(\kappa _{1},j)}=\{v_{a(1,1)}$, $v_{a(\kappa _{1},j,1)}$, $v_{a(\kappa _{1},j,2)}$, $\ldots$, $v_{a(\kappa _{2},j,k-1)}\}$ where $j=2, \ldots,  r_{1}$, $$B=\mathcal {F}_{(\mathcal {R}; \mathcal {S}; \mathcal {T})}-\sum^{r_{1}}_{j=2}e_{(\kappa _{1},j)}+\sum^{r_{1}}_{j=2}e^{'}_{(\kappa _{1},j)}$$ if $r_{1}\geq 2$. As Subcase 1.1 in Case 1, we get that $\rho(B)>\rho(\mathcal {F}_{(\mathcal {R}; \mathcal {S}; \mathcal {T})})$. If $r_{1}= 1$, then we let $B=\mathcal {F}_{(\mathcal {R}; \mathcal {S}; \mathcal {T})}$. Note that $B\cong \mathcal {F}_{2}$ now. Then using Lemma \ref{le01.01.01,00} gets $\rho(\mathcal {F})>\rho(\mathcal {F}_{(\mathcal {R}; \mathcal {S}; \mathcal {T})})$.

{\bf Subcase 1.2}  In $X$, $x_{v_{1}}\geq x_{v_{a(1,1)}}$. Now, if $s_{1}\geq 2$, then we let $e^{'}_{(\varepsilon_{1},j)}=\{v_{1}$, $v_{a(\varepsilon_{1},j,1)}$, $v_{a(\varepsilon _{1},j,2)}$, $\ldots$, $v_{a(\varepsilon _{1},j,k-1)}\}$ where $j=2, 3, \ldots, s_{1}$, $$B=\mathcal {F}_{(\mathcal {R}; \mathcal {S}; \mathcal {T})}-\sum^{s_{1}}_{j=2}e_{(\varepsilon_{1},j)}+\sum^{s_{1}}_{j=2}e^{'}_{(\varepsilon_{1},j)}$$ if $s_{1}\geq 2$. As Subcase 1.1 in Case 1, we get that $\rho(B)>\rho(\mathcal {F}_{(\mathcal {R}; \mathcal {S}; \mathcal {T})})$. If $s_{1}=1$, then we let $B=\mathcal {F}_{(\mathcal {R}; \mathcal {S}; \mathcal {T})}$. Note that $B\cong \mathcal {F}_{3}$ now. Then using Lemma \ref{le01.01.01,00} gets $\rho(\mathcal {F})>\rho(\mathcal {F}_{(\mathcal {R}; \mathcal {S}; \mathcal {T})})$.

{\bf Case 3} $r_{1}=0$, $r_{2}=0$. Note our previous supposition that if there exists some $s_{i}\geq 1$ or some $t_{i}\geq 1$ where $1\leq i\leq k-2$, then we suppose $s_{1}\geq 1$, $s_{i}=0$ for all $2\leq i\leq k-2$, and $t_{i}=0$ for all $1\leq i\leq k-2$. Thus $\mathcal {F}_{(\mathcal {R}; \mathcal {S}; \mathcal {T})}\cong \mathcal {F}_{1}$ now. Then using Lemma \ref{le01.01.01} gets $\rho(\mathcal {F})>\rho(F_{(\mathcal {R}; \mathcal {S}; \mathcal {T})})$.

From the above narrations, we conclude that $\rho(\mathcal {F})>\rho(F_{(\mathcal {R}; \mathcal {S}; \mathcal {T})})$. This completes the proof. \ \ \ \ \ $\Box$
\end{proof}

\begin{Proo}
Suppose $\mathcal{U}\in \Lambda$ satisfies that $\rho(\mathcal{U})=\max\{\rho(\mathcal{G})|\, \mathcal{G}\in\Lambda\}$.
Let $X$ be
the principal eigenvector of $\mathcal{U}$ and $x_{u}=\max\{x_{v}|\, v\in V(\mathcal{U})\}$. Denote by $\mathcal{C}=v_{1}e_{1}v_{2}e_{2}\cdots v_{q-1}e_{q-1}v_{q}e_{q}v_{1}$ the unique hypercycle in $\mathcal{U}$.

{\bf Assertion 1} $q=2$. Otherwise, suppose $q>2$.

{\bf Case 1} $q=3$. Now $\mathcal{C}=v_{1}e_{1}v_{2}e_{2}v_{3}e_{3}v_{1}$. Note that $\mathcal{U}$ is connected and $m(\mathcal{G})\geq 4$. Suppose $e_{1}=\{v_{1}$, $v_{a(1,1)}$, $v_{a(1,2)}$, $\ldots$, $v_{a(1,k-2)}$, $v_{2}\}$, $e_{3}=\{v_{3}$, $v_{a(3,1)}$, $v_{a(3,2)}$, $\ldots$, $v_{a(3,k-2)}$, $v_{1}\}$.

{\bf Subcase 1.1} There exists $v_{i}\in \{v_{1}$, $v_{2}$, $v_{3}\}$ such that $deg_{\mathcal{U}}(v_{i})\geq 3$. Without loss of generality, we suppose $deg_{\mathcal{U}}(v_{1})\geq 3$. If $x_{v_{2}}\geq x_{v_{3}}$, let $e^{'}_{3}=\{v_{2}$, $v_{a(3,1)}$, $v_{a(3,2)}$, $\ldots$, $v_{a(3,k-2)}$, $v_{1}\}$, $\mathcal{U}^{'}=\mathcal{U}-e_{3}+e^{'}_{3}$; if $x_{v_{2}}< x_{v_{3}}$, let $e^{'}_{1}=\{v_{1}$, $v_{a(1,1)}$, $v_{a(1,2)}$, $\ldots$, $v_{a(1,k-2)}$, $v_{3}\}$, $\mathcal{U}^{'}=\mathcal{U}-e_{1}+e^{'}_{1}$. Note that $\mathcal{U}^{'}\in \Lambda$. Using Lemma \ref{le03.02} gets that $\rho(\mathcal{U}^{'})>\rho(\mathcal{U})$, which contradicts the maximality of $\rho(\mathcal{U})$.

{\bf Subcase 1.2} There exists $v_{i}\in (V(\mathcal{C})\setminus\{v_{1}$, $v_{2}$, $v_{3}\})$ such that $deg_{\mathcal{U}}(v_{i})\geq 2$. Without loss of generality, then we suppose $deg_{\mathcal{U}}(v_{a(1,1)})\geq 2$. If $x_{v_{2}}\geq x_{v_{3}}$, let $e^{'}_{3}=\{v_{2}$, $v_{a(3,1)}$, $v_{a(3,2)}$, $\ldots$, $v_{a(3,k-2)}$, $v_{1}\}$, $\mathcal{U}^{'}=\mathcal{U}-e_{3}+e^{'}_{3}$; if $x_{v_{2}}< x_{v_{3}}$, let $e^{'}_{1}=\{v_{1}$, $v_{a(1,1)}$, $v_{a(1,2)}$, $\ldots$, $v_{a(1,k-2)}$, $v_{3}\}$, $\mathcal{U}^{'}=\mathcal{U}-e_{1}+e^{'}_{1}$. Similar to Subcase 1.1, we get $\mathcal{U}^{'}\in \Lambda$, and $\rho(\mathcal{U}^{'})>\rho(\mathcal{U})$, which contradicts the maximality of $\rho(\mathcal{U})$.

{\bf Case 2} $q\geq 4$. Suppose $e_{1}=\{v_{1}$, $v_{a(1,1)}$, $v_{a(1,2)}$, $\ldots$, $v_{a(1,k-2)}$, $v_{2}\}$, $e_{q}=\{v_{q}$, $v_{a(q,1)}$, $v_{a(q,2)}$, $\ldots$, $v_{a(q,k-2)}$, $v_{1}\}$. If $x_{v_{2}}\geq x_{v_{q}}$, let $e^{'}_{q}=\{v_{2}$, $v_{a(q,1)}$, $v_{a(q,2)}$, $\ldots$, $v_{a(q,k-2)}$, $v_{1}\}$, $\mathcal{U}^{'}=\mathcal{U}-e_{q}+e^{'}_{q}$; if $x_{v_{2}}< x_{v_{q}}$, let $e^{'}_{1}=\{v_{1}$, $v_{a(1,1)}$, $v_{a(1,2)}$, $\ldots$, $v_{a(1,k-2)}$, $v_{q}\}$, $\mathcal{U}^{'}=\mathcal{U}-e_{1}+e^{'}_{1}$. Similar to Subcase 1.1, we get $\mathcal{U}^{'}\in \Lambda$, and $\rho(\mathcal{U}^{'})>\rho(\mathcal{U})$, which contradicts the maximality of $\rho(\mathcal{U})$.

From the above two cases, we get that $q=2$. Then Assertion 1 follows as desired.

{\bf Assertion 2} There exists $\mu\in V(\mathcal{C})$ such that $x_{\mu}=x_{u}$. Otherwise, suppose $x_{v}<x_{u}$ for any $v\in V(\mathcal{C})$. Suppose $ue_{s_{1}}v_{s_{1}}e_{s_{2}}\cdots e_{s_{\varepsilon}}v_{s_{\varepsilon}}$ is the hyperpath from $u$ to $\mathcal{C}$ where $v_{s_{\varepsilon}}\in V(\mathcal{C})$. Note that $q=2$. Denote by $e_{1}=\{v_{1}$, $v_{a(1,1)}$, $v_{a(1,2)}$, $\ldots$, $v_{a(1,k-2)}$, $v_{2}\}$, $e_{2}=\{v_{1}$, $v_{a(2,1)}$, $v_{a(2,2)}$, $\ldots$, $v_{a(2,k-2)}$, $v_{2}\}$. Next, we consider two cases.

{\bf Case 1} $v_{s_{\varepsilon}}\in \{v_{1}, v_{2}\}$. Without loss of generality, suppose $v_{s_{\varepsilon}}=v_{1}$. Let $e^{'}_{2}=\{v_{2}$, $v_{a(1,1)}$, $v_{a(1,2)}$, $\ldots$, $v_{a(1,k-2)}$, $u\}$, $\mathcal{U}^{'}=\mathcal{U}-e_{2}+e^{'}_{2}$. Note that $\mathcal{U}^{'}\in \Lambda$ because $\mathcal{U}^{'}$ contains cycle $ue_{s_{1}}v_{s_{1}}e_{s_{2}}\cdots e_{s_{\varepsilon}}v_{1}e_{1}v_{2}e^{'}_{2}u$. As Assertion 1, it proved that $\rho(\mathcal{U}^{'})>\rho(\mathcal{U})$, which contradicts the maximality of $\rho(\mathcal{U})$.

{\bf Case 2} $v_{s_{\varepsilon}}\in (V(\mathcal{C})\setminus\{v_{1}, v_{2}\})$.  Without loss of generality, suppose $v_{s_{\varepsilon}}=v_{a(1,1)}$. Let $e^{'}_{2}=\{v_{2}$, $v_{a(2,1)}$, $v_{a(2,2)}$, $\ldots$, $v_{a(2,k-2)}$, $u\}$, $\mathcal{U}^{'}=\mathcal{U}-e_{2}+e^{'}_{2}$. Note that $\mathcal{U}^{'}\in \Lambda$ because $\mathcal{U}^{'}$ contains cycle $ue_{s_{1}}v_{s_{1}}e_{s_{2}}\cdots e_{s_{\varepsilon}}v_{a(1,1)}e_{1}v_{2}e^{'}_{2}u$. As Assertion 1, it proved that $\rho(\mathcal{U}^{'})>\rho(\mathcal{U})$, which contradicts the maximality of $\rho(\mathcal{U})$.

From the above two cases, then Assertion 2 follows as desired. By Assertion 2, we suppose $x_{\eta}=x_{u}$ where $\eta\in V(\mathcal{C})$. For any $v\in (V(\mathcal{U})\setminus V(\mathcal{C})$, let $dist_{\mathcal{U}}(v,\mathcal{C})$ denote the length of the shortest path from $v$ to $\mathcal{C}$. Similar to Claim 2 in the proof of  Theorem \ref{th01.01}, we get the following Assertion 3.

{\bf Assertion 3} $dist_{\mathcal{U}}(\eta, v)\leq 3$ and $dist_{\mathcal{U}}(v,\mathcal{C})= 1$ for any $v\in (V(\mathcal{U})\setminus V(\mathcal{C})$.

Now, from the above narrations, we get that $\mathcal{U}\cong \mathcal {F}_{(\mathcal {R}; \mathcal {S}; \mathcal {T})}$ for some $(\mathcal {R}; \mathcal {S}; \mathcal {T})$. Then using Lemmas \ref{le01.01.01}-\ref{le01.01.01,01} gets that $\mathcal{U}\cong \mathcal {F}$.
This completes the proof. \ \ \ \ \ $\Box$
\end{Proo}

\small {

}

\end{document}